\documentclass{article}
\usepackage{pgfplots}
\usepackage{fullpage}
\usepackage{amsmath}
\usepackage{amsthm}
\usepackage{amssymb}
\usepackage{url}
\usepackage{graphicx}
\usepackage{verbatim}
\usepackage{url}
\usepackage{graphicx}
\usepackage[caption=false]{subfig}
\usepackage{float}
\usepackage[ruled]{algorithm2e}
\newtheorem{theorem}{Theorem}

\theoremstyle{definition}

\newtheorem{remark}{Remark}

\usepackage{multirow}
\usepackage[utf8]{inputenc}
\begin{document}
\title{The image of the Veronese mapping is a determinantal variety; an elementary proof}
\author{Rahim Zaare-Nahandi\footnote{Emeritus Professor of Mathematics, University of Tehran, Tehran, Iran; rahimzn@ut.ac.ir}}
%\thanks{Emails: rahimzn@ut.ac.ir} 
%\keywords{Veronese imbedding, Determinantal  variety, Ideal of minors}
%\subjclass[2020]{14A10, 14A25}
\date{}
\maketitle

\begin{abstract}
\noindent
Let $\nu_d : \mathbb{P}^n \longrightarrow \mathbb{P}^N$ be the Veronese mapping of degree $d$ where $N = {n+d \choose n} -1$. By an elementary approach it is shown that $\nu_d$ is an isomorphism of $\mathbb{P}^n$ onto the projective variety defined by the $2$-minors of a matrix.
\end{abstract}

\section{Introduction}
The Veronese mapping, supposedly introduced by the Italian mathematician Giuseppe Veronese (1854-1917), is the map  $\nu_d : \mathbb{P}^n \longrightarrow \mathbb{P}^N$ given by 
$$[x_0 : x_1 : \cdots : x_n] \mapsto [x_0^d : x_0^{d-1}x_1 : x_0^{d-1}x_2 : \cdots : x_{n-1}x_n^{d-1} : x_n^d].$$ It is well-known that this map is an imbedding, i.e., its image is a projective variety, the so called {\it Veronese variety}, and the map is an isomorphism between $\mathbb{P}^n$ and the image variety. The Veronese imbedding is a basic tool in algebraic geometry (e.g., see 
the Indexes of \cite{Hart}, \cite{S} and \cite{Harr}). It allows one to reduce the study of some problems concerning  hypersurfaces of degree $d$ to the case of hyperplanes \cite[4.4, Example 2]{S}. 
%The image of $\nu_d$ is usually called the Veronese variety of degree $d$. 
The images of a fixed projective space under Veronese imbeddings of different degrees are isomorphic, but their homogenous coordinate rings are not. I.e., unlike the case of affine varieties, the homogenous coordinate ring is not a projective invariant.  
The Veronese imbedding has applications in other areas such as Coding Theory \cite{GP}.\\
 
%There are several known results in algebraic geometry, the so called ``folklores'', which require difficult or tediously long proofs where most readers do not bother reproving them. 
In most text books, when discussing the classical algebraic geometry, a rigorous proof of the fact that the image of the Veronese mapping is indeed a projective variety isomorphic to the initial projective space, is missing.  Some of these books leave the proof to the reader \cite[Exercises 2.12 and 3.4]{Hart} and some give sketches of a proof \cite[4.4, Example 2]{S}, \cite[5.1, Proposition]{SKKT} (see Remark 1). 
%Shafarevich  has given the main steps of the proof without much details; the defining equations are presented, it is hinted how certain open subsets cover the image of the Veronese mapping, and the inverse morphism is piecewise defined, whereas, checking most of these steps requires computations even for mappings of low degree.  
%In this regard, no criticism of these books is warranted as their purpose is to treat a wide range of topics and it would not make sense for all details to be supplied.\\
Nevertheless, no criticism of treatments of this result in these texts is warranted as their purpose is to treat a wide range of topics and it would not make sense for all details to be supplied.\\

The purpose of this note is to follow an elementary approach suitable to a curious student of a first course in algebraic geometry, and prove that the Veronese mapping is an isomorphism of $\mathbb{P}^n$ onto the projective variety defined by the $2$-minors of a $(n+1)\times {n+d-1 \choose n}$ matrix in the coordinates $[z_{i_0,\dots,i_n}: i_0+\dots +i_n = d]$ of 
$\mathbb{P}^N$. The matrix presentation helps in a systematic way, to provide a rigorous proof of the result. In the context of commutative algebra, in a recent paper by Abdolmaleki and Zaare-Nahandi \cite[Theorem 1.5 and Corollary 1.10]{AZ}, given an ideal $I$ of certain type generated by quadratic binomials, a matrix is introduced such that the set of its $2$-minors is a generating set of $I$. This paper has been the motivation of this note, since its treatment includes the defining ideal of the Veronese variety.  Primarily, in a completely different approach compared to \cite{AZ}, Pucci \cite[Corollaries 2.4 and 3.5]{P} have proved that, given the generators of the defining ideal of the Veronese variety, the rank-{1} determinantal ideal of certain catalecticant matrix provides the same generators. In the case of the Veronese mapping, the matrix obtained in \cite{AZ} is the same matrix  given in \cite{P}. 

%It needs mentioning that, by the term ``determinantal variety'' we mean a variety defined by $t$-minors of a matrix for some integer $t$. In the literature, this term has also been used for varieties defined by $t$-minors of a generic matrix, i.e., a matrix where all its entries are distinct indeterminates \cite{BV}.

\section{The construction and  the proof}

\noindent
{\bf Construction} \cite[Notation 1.3]{AZ}. Order the set of all monomials of degree $d$ in $x_0, \dots, x_n$ by the lexicographic order induced by  $x_0> x_1> \cdots >x_n$ . Let $L$ be the matrix where the entries of its $i$th row are all monomials of degree $d$ in $x_0, \dots, x_n$ with $x_i$ as a factor written decreasingly by the induced lexicographic order. Let $M$ be the matrix obtained from $L$ by replacing $x_0^{i_0} \dots x_n^{i_n}$ with $z_{i_0,\dots,i_n}$, the coordinates in $\mathbb{P}^N$. It is clear that both $L$ and $M$ are $(n+1)\times {n+d-1 \choose n}$ matrices.\\ 
%It is the same as the catalecticant matrix Cat$(1,d-1;n+1)$ in the notation of \cite{P}.\\

\noindent
{\bf Example.} Consider the imbedding 
$$\nu_3 : \mathbb{P}^2 \longrightarrow \mathbb{P}^9,$$  
$$[x_0 : x_1 : x_2] \mapsto [x_0^3: x_0^2x_1: x_0^2x_2: x_0x_1^2: x_0x_1x_2: x_0x_2^2: x_1^3: x_1^2x_2 : x_1x_2^2: x_2^3].$$
Then\\
\begin{center}
$L = \left(\begin{array}{cccccc}x_0^3 & x_0^2x_1 & x_0^2x_2 & x_0x_1^2 & x_0x_1x_2 & x_0x_2^2 \\x_0^2x_1 &x_0x_1^2& x_0x_1x_2 & x_1^3 & x_1^2x_2 & x_1x_2^2 \\ x_0^2x_2 &  x_0x_1x_2 &x_0x_2^2 & x_1^2x_2 & x_1x_2^2 & x_2^3\end{array}\right),$
\end{center}
and\\
\begin{center}
$M = \left(\begin{array}{cccccc}z_{3,0,0}& z_{2,1,0} & z_{2,0,1} & z_{1,2,0} & z_{1,1,1} & z_{1,0,2} \\z_{2,1,0} & z_{1,2,0} & z_{1,1,1} & z_{0,3,0} & z_{0,2,1} & z_{0,1,2} \\z_{2,0,1} &  z_{1,1,1} &z_{1,0,2} & z_{0,2,1} & z_{0,1,2} & z_{0,0,3}\end{array}\right).$
\end{center}

%\begin{thm}
\begin{theorem}
%\noindent
%{\bf Theorem 2.1.} 
Let $\nu_d : \mathbb{P}^n \longrightarrow \mathbb{P}^N$ be the Veronese mapping of degree $d$. Let $I_2(M)$ be the ideal generated by the set of $2$-minors of the matrix $M$ constructed above. Let $W = \mathbb{V}(I_2(M))\subset \mathbb{P}^N$ be the projective variety defined by $I_2(M)$. 
Then $W = {\rm Im}(\nu_d)$ and  the Veronese mapping is an imbedding of $\mathbb{P}^n$ onto $W$.
%\end{thm}
\end{theorem}

\noindent
\begin{proof}
By the construction, $L$ and $M$ may be presented by the following matrices, respectively:

\begin{equation}L=
\begin{pmatrix}
%$$\left(\begin{array}{ccccccccc}
x_0^d &x_0^{d-1}x_1&\cdots&x_0^{d-1}x_n&x_0^{d-2}x_1^2
 &\cdots&x_0x_n^{d-1}\\
x_0^{d-1}x_1 & x_0^{d-2}x_1^2&\cdots&x_0^{d-2}x_1x_n&x_0^{d-3}x_1^3
 &\cdots&x_1x_n^{d-1}\\
x_0^{d-1}x_2 &x_0^{d-2}x_1x_2& \cdots&x_0^{d-2}x_2x_n&x_0^{d-3}x_1^2x_2
 &\cdots&x_2x_n^{d-1}\\
\cdots&\cdots&\cdots&\cdots&\cdots
 &\cdots&\cdots\\
\cdots&\cdots& \cdots& \cdots&\cdots&\cdots&\cdots\\
x_0^{d-1}x_n& x_0^{d-2}x_1x_n&\cdots&x_0^{d-2}x_n^2&x_0^{d-3}x_1^2x_n
 &\cdots&x_n^d\\
%\end{array}\right),$$
\end{pmatrix},  
\end{equation}

\begin{equation}M=
\begin{pmatrix}
%$$\left(\begin{array}{ccccccccc}   
z_{d,0,\dots,0} &z_{d-1,1,\dots,0}&\cdots&z_{d-1,0,\dots,0,1}&z_{d-2,2,\dots,0}
&\cdots&z_{1,0,\dots,d-1}\\
z_{d-1,1,\dots,0} &z_{d-2,2,\dots,0}&\cdots&z_{d-2,1,0,\dots,1}&z_{d-3,3,\dots,0}
&\cdots&z_{0,1,\dots,d-1}\\
z_{d-1,0,1,\dots,0} &z_{d-2,1,1,\dots,0}&\cdots&z_{d-2,0,1,\dots,1}&z_{d-3,2,1\dots,0}
&\cdots&z_{0,0,1,\dots,d-1}\\
\cdots&\cdots&\cdots&\cdots&\cdots
&\cdots&\cdots\\
\cdots&\cdots&\cdots&\cdots&\cdots
&\cdots&\cdots\\
z_{d-1,0,\dots,1} &z_{d-2,1,0,\dots,1}&\cdots&z_{d-2,0,\dots,2}&z_{d-3,2,0\dots,i}
&\cdots&z_{0,0,\dots,d}\\
%\end{array}\right).$$
\end{pmatrix}. 
\end{equation}

\noindent
We proceed the proof in the following three steps:\\

\noindent
{\bf Step 1.}  
Let $W_i \subset W$, $i=0,1,\dots,n$, be the open subset of points where the coordinate $z_{0,\dots,0,d,0,\dots,0}$, corresponding to $x_i^d$, is nonzero. We show taht $W_0, \dots, W_n$ cover $W$. We prove that if $W_i$'s do not cover $W$, i.e., at some point $Q\in W$ all coordinates  $z_{0,\dots,0,d,0,\dots,0}$, corresponding to $x_i^d$, $i=0\dots, n$, are zero,   then all coordinates $[z_{i_0,\dots,i_n}: i_0+\dots +i_n = d]$ are zero, which would be a contradiction. 

Let's consider the $t$-th row of $M$, $t=0, \cdots, n-1$, corresponding to the $t$-th row of $L$ consisting of multiples of all monomials in $x_0, \cdots, x_n$ of degree $d-1$ by $x_t$. By hypothesis $z_{0,\dots, d, \dots, 0}$ corresponding to $x_t^d$ is zero. Assume that the vanishing of $2$-minors of $M$ implies that all entries on the $t$-th row of $M$ following $z_{0, \dots, d, \dots, 0}$ and preceding the entry $z_{j_0,\dots, j_t, \dots, j_n}$ corresponding to $x_0^{j_0}\dots x_t^{j_t}\dots x_n^{j_n}$, for some $j_0,\cdots,j_t, \cdots, j_n$, $j_t \ge 1$, are zero. Let $j_k$, $k>t$, be the last nonzero power in $x_0^{j_0}\dots x_t^{j_t}\dots x_n^{j_n}$. Thus, $x_0^{j_0}\dots x_t^{j_t}\dots x_n^{j_n}=x_0^{j_0}\dots x_t^{j_t}\dots x_k^{j_k}$. Consider the monomial $x_0^{j_0}\dots x_t^{j_t+1} \dots x_k^{j_k-1}$ obtained by reducing the power of $x_k$ by one and adding one to the power of $x_t$ in  $x_0^{j_0}\dots x_t^{j_t}\dots x_k^{j_k}$. Observe that  $x_0^{j_0}\dots x_t^{j_t+1} \dots x_k^{j_k-1}$ appears before $x_0^{j_0}\dots x_t^{j_t}\\ \dots x_k^{j_k}$ but not before $x_t^d$ on the $t$-th row of $L$ in the lexicographic order. Consider the $k$-th row of $L$, the row consisting of multiples of $x_k$. Recall that $t\le n-1$, therefore, the $k$-th row is below the $t$-th row of $L$. The entry on this $k$-th row and the column containing $x_0^{j_0}\dots x_t^{j_t+1} \dots x_k^{j_k-1}$ is  $x_0^{j_0}\dots x_t^{j_t}\dots x_k^{j_k}$.
Consider the $2$-minor of $M$
\begin{equation}m=
\begin{vmatrix} 
z_{j_0,\dots,j_t+1,\dots, j_k-1} &z_{j_0,\dots, j_t,\dots ,j_k} \\
z_{j_0,\dots, j_t,\dots ,j_k} & z_{j_0,\dots,j_t-1,\dots j_k+1} \\
\end{vmatrix}
\end{equation}  
corresponding to the $2$-minor of $L$
\begin{equation}
   \begin{vmatrix} 
x_0^{j_0}\dots x_t^{j_t+1}\dots x_k^{j_k-1} &x_0^{j_0}\dots x_t^{j_t}\dots x_k^{j_k}\\
x_0^{j_0}\dots x_t^{j_t}\dots x_k^{j_k}&x_0^{j_0}\dots x_t^{j_t-1}\dots x_k^{j_k+1} \\
\end{vmatrix}.
\end{equation}
By the lexicographic order,  $ z_{j_0,\dots,j_t+1,\dots, j_k-1} $ appears before $z_{j_0,\dots, j_t, j_k} $ and after  $z_{0, \dots, d, \dots, 0}$ on the $t$-th row of $M$. Thus, it is zero by the assumption. Hence, the vanishing of $m$ implies that $z_{j_0,\dots,j_n}=0$. Therefore, all entries on the $t$-row of $M$ starting from $z_{0, \dots, d, \dots, 0}$ are zero for all $t=0, \cdots, n-1$. Now, assume that the entries on the $i$-th row of $M$, for all $0\le i\le t-1$, preceding $z_{0, \dots, d, \dots, 0}$ corresponding to $x_i^d$ are also zero. The entries on the $t$-th row of $L$ before $x_t^d$, $t=1, \cdots, n$, are multiples of $x_i$ for some $0\le i\le t-1$. Thus, these entries appear on the previous rows of $L$. Therefore, the entries of $M$ corresponding to these entries of $L$ are already zero. It remain $z_{0, \dots, 0,d}$, which is zero by the assumption. Therefore, all coordinates $z_{i_0,\dots,i_n}$ are zero.\\

\noindent
{\bf Step 2.} We define maps $\varphi_i: W_i \longrightarrow \mathbb{P}^n$ for $i=0,\dots,n$ and show that these maps patch together to define a morphism $\varphi : W \longrightarrow \mathbb{P}^n$. Let $\varphi_i : W_i \longrightarrow \mathbb{P}^n$ be given by 
$$\varphi_i ([ z_{i_0,\dots,i_n}, i_0+\dots+i_n=d]) = {\rm the \ column \ of} \ M \  {\rm containing} \  z_{0,\dots, 0,d,0,\dots, 0},$$ 
where $z_{0,\dots, 0,d,0,\dots, 0}$ corresponds to $x_i^d$. We claim that, for all $i\not= j$, $\varphi_i$ and $\varphi_j$ agree on $W_i\cap W_j$. We do this for the columns containing  $z_{d,0,\dots,0}$ and $z_{0,d,0,\dots,0}$, the columns corresponding to $x_0^d$ and $x_1^d$, respectively. The procedure for other relevant columns is similar. In fact, the vanishing of the $2$-minors of the sub-matrix of $M$ consisting of the above two columns implies that these two columns are proportional, presenting the same point in $\mathbb{P}^n$. Therefore, the maps $\varphi_i$ match together to give a morphism $\varphi : W \longrightarrow \mathbb{P}^n$.\\

\noindent
{\bf Step 3.} We show that $\varphi$ is the inverse of $\nu_d$. This implies that  ${\rm Im}(\nu_d) = W$ is closed and $\nu_d$ in an imbedding of $\mathbb{P}^n$ onto ${\rm Im}(\nu_d)$. The equality $\varphi \circ \nu_d = id_{\mathbb{P}^n}$ is straightforward. Let $[x_0: \dots: x_n] \in  \mathbb{P}^n$ where $x_i\not= 0$ for some $i$. Then $\nu_d([x_0: \dots: x_n]) \in W_i$. Thus, $\varphi(\nu_d([x_0: \dots: x_n]))= \varphi_i(\nu_d([x_0: \dots: x_n])) = 
[x_0x_i^{d-1}:\dots:x_i^d:\dots: x_i^{d-1}x_n] = [x_0: \dots: x_n]$. To show $\nu_d\circ \varphi = id_W$, let $Q = [ z_{i_0,\dots,i_n}, i_0+\dots+i_n=d]\in W_i$. We assume that $i=0$. The argument for other values of  $i$ is similar. Thus, $$\nu_d\circ \varphi (Q) = \nu_d\circ \varphi_0 (Q) = \nu_d([z_{d,0,\dots,0}: z_{d-1,1,0,\dots,0}: \cdots : z_{d-1,0,\dots,0,1}]) = $$ $$=[z_{d,0,\dots,0}^d: z_{d,0,\dots,0}^{d-1}z_{d-1,1,0,\dots,0}: \cdots :z_{d-1,0,\dots,0,1}^d].$$
Given any coordinate $z_{d,0,\dots,0}^{i_0} z_{d-1,1,0,\dots,0}^{i_1} \cdots z_{d-1,0,\dots,0,1}^{i_n}$ of $\nu_d\circ \varphi_0 (Q)$, using some $2$-minors of $M$ where three entries of each minor belong to the first row and the first column of $M$, we show that $$z_{d,0,\dots,0}^{i_0} z_{d-1,1,0,\dots,0}^{i_1} \cdots z_{d-1,0,\dots,0,1}^{i_n} = z_{d,0,\dots,0}^{d-1}z_{i_0, \dots, i_n}.$$ This would mean that $\nu_d\circ \varphi = id_W$. The claim is already valid for $i_0=d, d-1$. Thus we may assume that $i_0 < d-1$. Let $i_s$ be the last nonzero integer in $(i_1, \dots, i_n)$. Let $z_{d-1,0,\dots,0,1,0\dots,0}$ correspond to $x_0^{d-1}x_s$. Let $i_s\ge 2$. Using the minor of $M$
\begin{equation}
   \begin{vmatrix} 
 z_{d,0,\dots,0} &z_{d-1,0,\dots,0,1,0\dots,0} \\
z_{d-1,0,\dots,0,1,0\dots,0} & z_{d-2,0,\dots,0,2,0,\dots,0} \\
\end{vmatrix},
\end{equation}
corresponding to the minor
\begin{equation}
   \begin{vmatrix} 
x_0^d &x_0^{d-1}x_s\\
x_0^{d-1}x_s&x_0^{d-2}x_s^2 \\
\end{vmatrix}
\end{equation}
of $L$, it follows that  
$$z_{d-1,0,\dots,0,1,0\dots,0}^{i_s} =  z_{d,0,\dots,0}z_{d-1,0,\dots,0,1,0,\dots,0}^{i_s-2}z_{d-2,0,\dots,0,2,0,\dots,0}.$$
If $i_s-2 \ge 1$, then using the $2$-minor
\begin{equation}
   \begin{vmatrix} 
 z_{d,0,\dots,0} &z_{d-2,0,\dots,0,2,0\dots,0} \\
z_{d-1,0,\dots,0,1,0\dots,0} & z_{d-3,0,\dots,0,3,0,\dots,0} \\
\end{vmatrix},
\end{equation}
corresponding to the minor
\begin{equation}
   \begin{vmatrix} 
x_0^d &x_0^{d-2}x_s^2\\
x_0^{d-1}x_s&x_0^{d-3}x_s^3 \\
\end{vmatrix}
\end{equation}
of $L$, it follows that 
$$z_{d-1,0,\dots,0,1,0,\dots,0}z_{d-2,0,\dots,0,2,0,\dots,0}=  z_{d,0,\dots,0}z_{d-3,0,\dots,0,3,0,\dots,0}.$$
Therefore,
$$z_{d-1,0,\dots,0,1,0,\dots,0}^{i_s}=  z_{d,0,\dots,0}^2z_{d-1,0,\dots,0,1,0,\dots,0}^{i_s-3}z_{d-3,0,\dots,0,3,0,\dots,0}.$$
Finally, continuing similar procedure, since $i_s\le d$, using the $2$-minor
\begin{equation}
   \begin{vmatrix} 
 z_{d,0,\dots,0} &z_{d-i_s+1,0,\dots,0,i_s-1,0\dots,0} \\
z_{d-1,0,\dots,0,1,0\dots,0} & z_{d-i_s,0,\dots,0,i_s,0,\dots,0} \\
\end{vmatrix},
\end{equation}
corresponding to the minor
\begin{equation}
   \begin{vmatrix} 
x_0^d &x_0^{d-i_s+1}x_s^{i_s}\\
x_0^{d-1}x_s&x_0^{d-i_s}x_s^{i_s} \\
\end{vmatrix}
\end{equation}
of $L$, it follows that 
$$z_{d-1,0,\dots,0,1,0,\dots,0}z_{d-i_s+1,0,\dots,0,i_s-1,0,\dots,0}=  z_{d,0,\dots,0}z_{d-i_s,0,\dots,0,i_s,0,\dots,0}.$$
Therefore,
$$z_{d-1,0,\dots,0,1,0\dots,0}^{i_s} =  z_{d,0,\dots,0}^{i_s-1}z_{d-i_s,0,\dots,0,i_s,0,\dots,0}.$$
Now, let $i_r$ be the nonzero integer in $(i_1, \dots, i_n)$ right before $i_s$. Using the minor of $M$
\begin{equation}
\begin{vmatrix} 
z_{d,0,\dots,0} &z_{d-i_s,0,\dots,0,i_s,0\dots,0} \\
z_{d-1,0,\dots,0,1,0\dots,0} & z_{d-i_s-1,0,\dots,0,1,0,\dots,0,i_s,0,\dots, 0} \\
\end{vmatrix},
\end{equation}
corresponding to the minor
\begin{equation}
\begin{vmatrix} 
x_0^d &x_0^{d-i_s}x_s^{i_s}\\
x_0^{d-1}x_r&x_0^{d-i_s-1}x_rx_s^{i_s} \\
\end{vmatrix}
\end{equation}
of $L$ we get 
$$z_{d-1,0,\dots,0,1,0\dots,0}z_{d-i_s,0,\dots,0,i_s,0\dots,0}= z_{d,0,\dots,0}z_{d-i_s-1,0,\dots,0,1,0,\dots,0,i_s,0,\dots,0}.$$
Therefore,
$$z_{d-1,0,\dots,0,1,0\dots,0}^{i_r}z_{d-i_s,0,\dots,0,i_s,0\dots,0}= z_{d,0,\dots,0}z_{d-1,0,\dots,0,1,0\dots,0}^{i_r-1}z_{d-i_s-1,0,\dots,0,1,0,\dots,0,i_s,0,\dots,0}.$$
Continuing this procedure similar to the case of $i_s$, we reach
$$z_{d-1,0,\dots,0,1,0,\dots,0}^{i_r}z_{d-1,0,\dots,0,1,0,\dots,0}^{i_s}= z_{d,0,\dots,0}^{i_r+i_s-1}z_{d-i_r-i_s,0,\dots,0,i_r,0,\dots,0,i_s,0,\dots,0}.$$
Further procedure on the other nonzero powers of $z_{d-1,0,\dots,0,1,0,\dots,0}$'s corresponding to  $x_0^{d-1}x_t$ with $t<r$, it follows that
$$z_{d,0,\dots,0}^{i_0} z_{d-1,1,0,\dots,0}^{i_1} \cdots z_{d-1,0,\dots,0,1}^{i_n} = z_{d,0,\dots,0}^{d-1}z_{i_0, \dots, i_n}.$$
Observe that since the $(1,1)$-entry of every minor employed in the reduction procedure is $z_{d,0,\dots,0}$, three entries of each minor belong to the first row and the first column of $M$.\\

If $Q = [ z_{i_0,\dots,i_n}, i_0+\dots+i_n=d]\in W_i$, for $i\not= 0$, using some $2$-minors of $M$ where three entries of each minor belong to the $i$-th row and the $i$-th column of $M$, similar to the case $i=0$, it follows that $$z_{d,0,\dots,0}^{i_0} z_{d-1,1,0,\dots,0}^{i_1} \cdots z_{d-1,0,\dots,0,1}^{i_n} = z_{0,\dots,0,d,0\dots, 0}^{d-1}z_{i_0, \dots, i_n},$$
as required. 
\end{proof}

\begin{remark} 
%\noindent
%{\bf Remark 2.2.} 
The proof of Theorem 1 in \cite[5.1, Proposition]{SKKT} may need some explanation. It begins with defining the inverse map of $\nu_d$. It is important to mention that the image of $\nu_d$ is closed in the ambient projective space, namely, it is a projective variety. To show this, one way is to prove that the image of $\nu_d$ is given by certain equations.  Such equations are given in the book, but after the proof of the Proposition. This property is not valid for maps of affine varieties. As an example, the image of the hyperbola $xy-1=0$ under the projection to its $x$-coordinate is not closed in the affine line. It is in fact a fundamental property of maps of projective varieties that the image of a projective variety under a projective map is closed. This is a consequence of the so called ``Fundamental Theorem of Elimination Theory'' (e.g., see \cite[5.2, Theorem 2]{S} or \cite[Theorem 14.1 and Corollary 14.2]{E}). The proof in  \cite[5.1, Proposition]{SKKT} ends with stating that the composition map $W \longrightarrow \mathbb{P}^N \longrightarrow W$ is the identity map. Again, to show this, one needs to use the equations proposed for the image of $\nu_d$.\\
\end{remark}

\noindent
{\bf Acknowledgements}\\
The author would like to thank his brother Rashid Zaare-Nahandi for bringing attention to his joint paper \cite{AZ} which motivated this note. He thanks Hassan Haghighi for reading the first manuscript of the paper and for some useful comments.

\end{document}